 \documentclass[12pt,a4paper,leqno]{article}
 \usepackage{amsmath}
\usepackage{amssymb}
\usepackage{amsfonts}
 \usepackage{amsthm}
\usepackage[english, french]{babel}
\usepackage[applemac]{inputenc}
\usepackage[all]{xy}
\usepackage{mathrsfs}
\numberwithin{equation}{section}
\renewcommand{\phi}{\varphi}

\DeclareMathOperator{\Id}{Id}
\DeclareMathOperator{\trace}{trace}

\def \R {\mathbb R}

\begin{document}

\newtheorem{theorem}{Theorem}[section]
\newtheorem{lemma}[theorem]{Lemma}
\newtheorem{proposition}[theorem]{Proposition}
\newtheorem{Cor}[theorem]{Corollary}

\def \cal {\mathcal}

\theoremstyle{definition}
\newtheorem{definition}[theorem]{Definition}
\newtheorem{example}[theorem]{Example}
\newtheorem{exer}[section]{Exercise}
\newtheorem{conj}{Conjecture}

\theoremstyle{remark}
\newtheorem{remark}[theorem]{Remark}
\newtheorem{remarks}[theorem]{Remarks}
\newtheorem{question}{Question}

%
\selectlanguage{english}
\title{\textbf {On the uniqueness of generating Hamiltonian for continuous limits of Hamiltonians flows}}
\author{  \\{ \sc Claude Viterbo}
\\
 Centre de Math\'ematiques Laurent Schwartz \\ UMR 7640 du CNRS
\\École Polytechnique - 91128 Palaiseau, France
\\ {\tt viterbo@math.polytechnique.fr   } }
\maketitle

\begin{abstract}
We show that if a sequence of Hamiltonian flows has a $C^0$ limit, and if the 
generating Hamiltonians of the sequence have a limit, this limit is uniquely determined 
by the limiting $C^0$ flow. 
This answers a question by Y.G. Oh in \cite {Oh-continuous}.
\end{abstract}
\section{Introduction}

Let $H_k(t,x)$ be a sequence of Hamiltonians with flow $\phi_k^t$ on a
symplectic manifold $(M,\omega)$, such that $\lim_k H_k(t,x)=H(t,x)$ and
$\lim_k \phi^t_k=\phi^t$ where limits are intended as $C^0$ limits.
Can we say that $\phi$ is the flow of $H$ ? In the case where $H$
has a flow (e.g. $H$ is $C^{1,1}$) this has been  proved in  \cite{Viterbo-STAGGF}, and one could alternatively use the methods  of
\cite{Hofer-displ}.  However if $H$ is only $C^0$, it is not easy to make
sense of this question, since the flow of a $C^0$ hamiltonian is
not defined.

This question is not as artificial as the reader may think, and
 has apparently been asked by Y.G. Oh in the framework of
$C^0$-Hamiltonians (\cite{Oh-continuous}, this  seems to be related to Question 3.11 or 3.20). It is sufficient to solve this question for the case where $H$ is continuous and $\phi=Id$. Do we necessarily
have $H=0$ in this case ?  

Of course we could have $H(t,x)=h(t)$ and we exclude this case by normalizing the Hamiltonians, in the compact case by imposing the condition $$\int_M H_k(t,x) \omega^n =0$$  and in the non-compact case by assuming they have  compact support. 

Our aim in this short note is to give a positive answer to Oh's
question. We denote by $C^{1,1}$ the set of differentiable functions with Lipschitz derivative.

\begin{theorem}
Let $H_n(t,z)$ be a sequence of $C^{1,1}$ Hamiltonians on
$(M,\omega)$, normalized as above, and such that $H_n$ converges in the $C^0$ topology to
some continuous function $H(t,z)$. Let $\phi_n^t$ be the flow of
$H_n$. Then if $\phi_n^t$ converges to $\Id$ in the $C^0$
topology, we have $H=0$.
 \end{theorem}
 
  \begin{remark} 
a)  Throughout the paper, by $C^0$ convergence of $\phi_n^t$ to $\phi^t$, we always mean $C^0$  convergence uniform in $t$, for $t$ in $[0,1]$. In other words,
  
  $$ \forall \varepsilon >0\; \exists N_0\; \forall t \in [0,1] \; \forall  N>N_0,\;  \Vert   \phi_n^t- \phi^t \Vert_{C^0} \leq \varepsilon $$
  
  b) Note that if we do not assume that $H_n$ converges, the theorem does not hold. Indeed, consider a non-zero Hamiltonian $H_0$ supported in the unit ball. Then the sequence $nH_0(nz)$ does not converge, but the time one flow $C^0$, does converge to the identity.  

c) According to Y.G. Oh, one can adapt the proof of the theorem  to the case where convergence is in Hofer norm, i.e. the norm given by $$\Vert H \Vert = \int_0^1 \left [ \max_{x\in M}H(t,x)-\min_{x\in M} H(t,x) \right ] dt $$ 
 \end{remark} 

As a Corollary we get
 
 \begin{Cor} 
 Let $\phi_n^t, \psi_n^t$ be sequences of Hamiltonian flows associated to $H_n(t,x),K_n(t,x)$. 
 Assume $$\lim_n H_n=H, \lim_n K_n=K$$ and $$\lim_n\phi_n^t= \lim_n \psi_n^t = \rho^t$$
 where all limits are intended as $C^0$ limits. Then 
 $$H=K$$
 \end{Cor} 
   
   \begin{proof}[Proof of Corollary, assuming the theorem] 
   Indeed, $(\phi_n^t)^{-1} \circ (\psi_n^t)$ $C^0$ converges to the identity, and  is generated by $(K_n-H_n)(t, \phi_n^t(z))$, hence $C^0$ converges to $(H-K)(t,\rho^t(z))$. Thus, 	according to the theorem we have $(H-K)(t,\rho^t(z))=0$ hence $H=K$. 
   
   \end{proof} 
  
 We thank Albert Fathi for drawing our attention to this problem, and Y.G. Oh for raising the question and for some useful comments. 
 
\section{On $C^0$-limits of Lagrangian submanifolds}



Consider  the problem of a topological
submanifold (i.e. $C^0$) $L$ in $T^*N$, that would be a $C^0$
limit of  $C^1$ Lagrangian submanifolds. According to
\cite{Laudenbach-Sikorav-IMRN}, if $L$ is $C^1$, it is necessarily
Lagrangian. When $L$ is  $C^0$, the meaning of  ``Lagrangian'' is
unclear. However if $L$ is a graph in the cotangent bundle, $L=\{
(x, p(x)) \mid x \in N\}$, requiring that  $p(x)dx$ is closed
makes sense even if $p\in C^0$. Indeed we may interpretate this as
meaning that $p(x)dx$ is closed in the distribution sense, as
suggested by Michael Herman (\cite{Herman} definition 8.13 page
60). In our case we wish to prove that if $L$ is a non Lagrangian
$C^0$ graph, we may not approximate it by Lagrangian submanifolds.
 Our
crucial assertion is
\begin{proposition}
 Let $N$ be a closed manifold  that is the total space of an $S^1$
fibration.
 Let  $p$ be a continuous section of $T^*N$ which, considered as a
one-form,
 is not closed in the sense of distributions. Then, there exists
$f\in C^{\infty}$ such that
 $p(x)-df(x)$ does not vanish on $N$.
 \end{proposition}

   \begin{Cor}\label{Corollary}
Let  $L$ be the graph of a one form $p$ and assume that any
neighbourhood $U$ of $L$ contains a smooth exact Lagrangian submanifold of
$T^*M$. Then $p$ is closed in the
sense of distributions.

In particular if  $L_n$ is a sequence of exact Lagrangians submanifolds in
$T^*M$,  $L$ is the $C^0$ graph of a
   one form $p$, and $L_n$ converges $C^0$ to $L$, then $p$ is closed.
   \end{Cor}
   
   This corollary is proved at the end of this section. 
 \begin{remark} \begin{enumerate}
    \item Note that once we have a $C^1$ function $f$ such that  $p(x)-df(x)$ does not vanish, the same holds if we replace $f$ by any $g$ close to $f$. Thus, we can replace $f$ by a smooth
approximation.
   In the sequel we shall thus not bother about the smoothness of the solution.
    \item Note that a smooth fibration with fiber $S^1$, that is a
principal bundle with group $Diff(S^1)$
    is equivalent to a principal bundle with fiber the group $S^1$,
    since the inclusion of $S^1$ into $Diff(S^1)$ induces an
    isomorphism of homotopy groups.
    
    \item  Our corollary can be compared with Theorem 2 in \cite{Sik3}. Even though it is stated
    there with stronger assumptions, Sikorav's proof yields  our corollary when $p$ is smooth.
    It can probably be adapted to  the case where $p$ is only $C^0$ using an analog of our
    lemma 2.5.  
    
Note that in  theorem 2 of \cite{Laudenbach-Sikorav-IMRN} the authors consider a sequence 
$\phi_n:V \longrightarrow T^*M$ of maps converging $C^0$  to $\phi_\infty$. They 
prove that if the $L_n=\phi_n(V)$ are Lagrangians embeddings, and $\phi_\infty$ 
is smooth, then $\phi_\infty(V)$ is Lagrangian, provided $\pi_2(T^*M,L_n)=0$. 
In particular, up to a symplectomorphism, $L_n$ may be assumed to be exact. 
In this case however, the smoothness of $\phi$ seems to be crucial, and moroever
convergence here is meant in a stronger sense (convergence of embeddings rather
than Hausdorff). 
    
Our Corollary trades the exactness requirement against a weaker assumption on convergence. 
This  exactness of $L_n$ is crucial as can be seen from the following example. 
    
Any submanifold can be approximated in the Hausdorff
topology by a (non-exact)  Lagrangian one.  Indeed, given $V$, we may
approximate it by a union of small Lagrangian tori, each being
contained in a Darboux chart near $V$. On the union of such tori,
we may perform a Polterovich surgery, in order to obtain a
connected Lagrangian submanifold. 

Note that quite obviously the map
$H_1(V)\to H_1(N)$ is not injective in this case: each torus
produces a lot of $1$-cycles of $V$, which go to zero in $N$.
This non-injectivity of the map $H_1(V)\to H_1(N)$ is crucial in this  counterexample:
if the map was injective, we could translate $\L_n$ by a closed form to make it exact,
and our Corollary would apply.
 \end{enumerate} 
 \end{remark}
 
 \bigskip  The proof of the proposition will require some lemmata. \newline
 Under the  assumptions of the proposition, let  $V$ be the base of the circle
fibration, and $y:N \longrightarrow V$ be the projection. For $y$ in
  a  domain of trivialization of the fibration, we consider
coordinates $(\theta,y)$ where $\theta\in S^1$.
  Note that we may always assume to be given an invariant measure (by
the circle action) on $N$. We also set
 $p(\theta,y)=(\pi (\theta,y), r(\theta,y))$ where $\theta\in S^1,
y\in V$, $\pi(\theta,y)\in {\mathbb R}$
 and $r$ is a section of $T^*V$ parametrized by $\theta\in S^1$.

 \begin{lemma}
 Assume  for some $y \in V$ we have
 $$P(y)=\int_{S^1} \pi(\theta,y) d\theta\neq 0 $$
 Then there exists $f\in C^1(S^1\times V,\R)$ such that
$p(\theta,y)-df(\theta,y)$ does not vanish.
 The same holds for $N$ the total space of a circle fibration..
 \end{lemma}
 \begin{proof}
Let $p(\theta,y)=(\pi(\theta,y), r(\theta,y))$ be a one form on $S^1\times
V$ that is a section of $T^*(S^1\times V)$. We look for a function
$f(\theta,y)$ such that $p(\theta,y)-df(\theta,y)$ never vanishes.
 Let us first try to solve $$ \pi(\theta,y)- \frac{ \partial f}{
\partial \theta}(\theta,y) = \varepsilon (\theta,y)$$
 Then this is solvable if and only if $$\int_{S^1} (\pi(\theta,y)-
\varepsilon (\theta,y))d\theta=0$$
 and thus, denoting by $P(y)=\int_{S^1}\pi(\theta,y)d\theta$, we can
choose $ \varepsilon$ non vanishing outside a neighbourhood of the set
$Z=\{(\theta,y) \mid P(y)=0\}$ (e.g. take $ \varepsilon (\theta,y)=P(y)$). In general we cannot choose $
\varepsilon$ to be non-zero in such a neighbourhood.

 Note that $f$ is well-defined up to a function of $y$. Also, we can assume $f$ to be smooth, provided we took care to choose $ \varepsilon (\theta,y)-\pi(\theta,y)$ smooth. 
 

 Now we need to find $h\in C^1(V,\R)$ such that $p(\theta,y)-df(\theta,y)-dh(y)$
does not vanish on $Z$. But if the projection of $Z$ on $V$, $U$, is
not all of $V$,  we can find a function $h$ on $V$ with no critical
point in $U$. Multiplying $h$ by a large constant, we may assume $dh$
to be arbitrarily large. Then $p(\theta,y)-df(\theta,y)-dh(y)$ will
not vanish for all $(\theta,y)$ with $y\in U$, and thus for
$(\theta,y)\in Z$.

 \end{proof}

\begin{lemma}
Assume $p$ is $C^0$ on $N$, and consider a smooth circle fibration of
$N$. Assume for any curve $\gamma$, $C^\infty$  close to a fiber of the fibration, we have
$$\int_{S^1} p(\gamma(t))\dot\gamma (t) dt=0$$
Then $p$ is closed in the sense of distributions.
\end{lemma}

 \begin{proof}

 We shall take local coordinates $(\theta,y)$ in the neighbourhood of a fiber, and let
$\eta(\theta,y)$ be a smooth vector  field on $N$.

 Let $\alpha$ be a continuous one-form. We wish to compute the integral of
$\alpha$ over the curve
 $t \to ( t, y+  \varepsilon \eta(t , y) )$. This will
be

 $$\int_{S^1}\alpha( \theta, y+ \varepsilon \eta(\theta , y ))
 (1, \varepsilon  \frac{\partial }{ \partial \theta }\eta(\theta , y
)) d \theta
 $$

 Writing $\alpha = \alpha_\theta d \theta + \alpha_ydy$ we can
rewrite the above as

 $$
 \int_{S^1}\alpha_ \theta ( \theta, y+ \varepsilon \eta(\theta , y )
) +  \varepsilon
 \alpha_ y (\theta, y+ \varepsilon \eta(\theta , y ) )
\frac{\partial }{\partial \theta }\eta(\theta , y ) d \theta
 $$

 Now averaging this over  $y \in V$, and differentiating with respect to $
\varepsilon $, we get

 \begin{gather*}
 \frac{\partial }{\partial \varepsilon }  \int_V \int_{S^1}\alpha_
\theta (\theta, y+ \varepsilon \eta(\theta , y ) ) +  \varepsilon
\alpha_ y (\theta, y+ \varepsilon \eta(\theta , y )  )
\frac{\partial }{\partial \theta }\eta(\theta , y ) d \theta dy = \\
-  \int_V \int_{S^1}\bigg [ \alpha_ \theta (\theta, y+
\varepsilon \eta(\theta , y ) ) \nabla_y\cdot \eta(\theta ,y ) +   \alpha_
y (\theta, y+ \varepsilon \eta(\theta , y ))  \frac{\partial
}{\partial \theta }\eta(\theta , y )+   \\ \varepsilon   \alpha_y  (\theta, y+ \varepsilon \eta(\theta , y )) \wedge d_y\eta (\theta,y)\bigg ] d\theta dy
 \end{gather*}
This is obtained by applying the change of variable $(\theta',y')=(\theta, y+ \varepsilon \eta(\theta,y))$ to the 

$$\int_{V\times S^1} f(\theta, y+ \varepsilon \eta(\theta,y))  d\theta dy$$

to get $$\int_{V\times S^1} f(\theta,y') dy$$ where $$dy'=dy + \varepsilon \eta(\theta,y') + o( \varepsilon )$$ so that $$dy=dy'- \varepsilon \eta(\theta, y') + o ( \varepsilon )$$
and $$\det (dy')= \det (dy) - \varepsilon \nabla_{y} \cdot \eta (\theta,y)$$
remembering that $ \trace d \eta = \nabla \cdot \eta$
$$\int_{V\times S^1} f(\theta,y') dy=\int_{V\times S^1} f(\theta,y') dy'- \varepsilon \int_{V\times S^1} f(\theta,y') \frac{\partial}{\partial y'}\eta (\theta,y') dy'  $$ and denoting by $\nabla_y$ the nabla operator with respect to the $y$ variables,  where all derivatives should be understood in the distributional
sense. Now taking  the above for $ \varepsilon =0$, we get

  \begin{gather*}
  \int_V \int_{S^1} \left [  (\alpha_\theta
\nabla_y\cdot \eta)(\theta ,y ) +   \alpha_ y(\theta , y )
\frac{\partial }{\partial \theta }\eta(\theta , y ) \right ] d \theta
dy
   \end{gather*}

 Integrating by parts we get

   \begin{gather*}
  \int_V \int_{S^1} \left [ (\nabla_y  \alpha_
\theta) (\theta , y ) \eta(\theta , y ) -  (\frac{\partial }{\partial
\theta }    \alpha_ y )(\theta , y )   \eta(\theta , y) \right ] d
\theta dy = \\
    \int_V \int_{S^1} \left [( \nabla_y  \alpha_
\theta) (\theta , y ) -  \frac{\partial }{\partial \theta }
\alpha_y (\theta , y)    \right ]  \eta(\theta , y )d \theta dy
   \end{gather*}
%
 
The last line equals the integration of $d\alpha$ against the
bivector  $\frac{\partial}{\partial \theta}\wedge (0, \eta) $. As this vanishes for all $\eta$, means that
${\imath}_{ \frac{\partial }{\partial \theta }}d \alpha $ vanishes as
a distribution (or current).

We thus proved that if for all $\eta$ the integration of $\alpha$
over the loop $t\to (t, y+  \varepsilon \eta(t , y))$ has vanishing derivative, then
${\imath}_{\frac{\partial }{\partial \theta }}d \alpha $ is identically zero.
Now if we slightly modify our fibration, and apply the same
argument, we get that ${\imath }_{Z}d \alpha =0$ for any  vector
field $Z$ tangent to
the  fiber
of a circle  fibration of $N$, close to the given one. The next lemma allows us to conclude the proof.

 \begin{lemma} 
 Assume $\alpha$ is a continuous form such that for any  vector field $Z$, tangent to a circle fibration of $N$ and close to $Z_0$, we have $i_Zd\alpha=0$. Then $d\alpha=0$ as a distribution.
  \end{lemma} 
   \begin{proof} 
   Indeed, it is enough to show that our assumption implies that $i_Zd\alpha$ vanishes for all vector fields $Z$. 
   
First of all, the problem is local: using a partition of unity, it is enough to show that 
$i_Zd\alpha =0$ holds for any $Z$ supported in a small set, tangent to a fibration close to $Z_0$. 

Now since $Z_0$ does not vanish, any vector field $C^1$ close to $Z_0$ has a flow box near $z_0$, hence a small diffeomorphism makes it tangent to $Z_0$. Thus locally, the set of $Z$ such that $i_Zd\alpha=0$ is open in the $C^\infty$ topology, and thus, by considering $i_{Z-Z_0} $, any small vector field supported  in the neighbourhood of $z_0$ satisfies $i_Zd\alpha=0$.    \end{proof} 
  \end{proof}
  \begin{proof}[Proof of the proposition]
 According to the second lemma, if $p$ is not closed, using a
vector field, we may smoothly perturb the fibration, $\pi$ so that
one of the fibers satisfies $\int_{\pi^{-1}(y)}p \neq 0$. Then, using
this new fibration and the first lemma, we see that there is a
function $f$ such that $p(x)-df(x)$ does not vanish.

  \end{proof}

   \begin{proof} [Proof of Corollary, following
\cite{Laudenbach-Sikorav-IMRN}]
First of all  if $L_n$ converges to $L$, then  $L_n\times 0_{S^1}
\subset T^*(N\times S^1)$
converges to $L\times 0_{S^1}$, and this will be the graph of $p$,
considered as a one-form
on $N\times S^1$. Now if $p$ is closed on $N$, its extension to
$N\times S^1$ is also closed, since
$$\int_{N\times S^1} \left [\frac{\partial}{\partial x_i} p_j(x)-
\frac{\partial}{\partial x_j}p_i(x)
\right ] \phi(x,\theta) dx d\theta $$ defined as

$$- \int_{N\times S^1} \left [p_j(x)\frac{\partial}{\partial x_i}
\phi(x,\theta)- p_i(x) \frac{\partial}{\partial x_j}\phi(x,\theta)
\right ]
 dx d\theta $$
 is equal to
$$ \int_{N\times S^1} \left [p_j(x)\frac{\partial}{\partial x_i}
\bar\phi(x)- p_i(x) \frac{\partial}{\partial x_j}\bar\phi(x)\right ]
dx d\theta$$

where we set $\bar\phi (x)=\int_{S^1}\phi(x,\theta)d\theta$, so that
$p$ is closed (in the sense of distributions)  as a one form on $N$
if and only if it is closed
(in the sense of distributions) as a one form on $N\times S^1$.

 According to the above lemma, we see that we may, using a
Hamiltonian symplectomorphism,
 send  $L$ away from the zero section (by $(x,p)\to (x, p-df(x))$)
and thus
 any Lagrangian submanifold $L_n$ in a neighbourhood of $L$ will also
be sent to
 $T^*N\setminus 0_N$ and thus may be disjoined from itself by a
small Hamiltonian isotopy,
 since $(x,p)\to (x,\lambda p)$ is conformal, and thus induces a
Hamiltonian isotopy on exact Lagrangians.
 But this is impossible according to Gromov's theorem (\cite{Gromov}
p. 330). \end{proof}

\section{Proof of the  theorem}
 
 \begin{lemma} 
Let  $K(t,z)$ be a Hamiltonian in $T^*N$, with flow $\psi^t$. Then  the embedding

 \begin{gather*} \Psi: [0,1]\times N \longrightarrow T^*([0,1]\times N)\\
(t,z)  \longrightarrow (t,-K(t,\psi^t(z)),\psi^t(z)) \end{gather*} 

is exact Lagrangian.  If moreover $\psi^{t}(z)$ and $K(t,z)$ are $1$-periodic in $t$, the Lagrangian submanifold of $T^*(S^1\times N)$ thus obtained is also exact, after maybe changing $K$ by a constant. 
\end{lemma} 

 \begin{proof} 
Indeed, if $\lambda $ is the Liouville form, and denoting by $d$ the differential with respect to $x$, while $D$ is the differential with respect to both $t$ and $x$, 

$$\Psi^*(\lambda+hdt)= (\psi^t)^*(\lambda)+(\lambda(\frac{d}{dt}\psi^t(z)))-K(t, \psi^t(z)) dt 
$$
$$= (\psi^t)^*(\lambda)+ (\psi^t)^*[i_{X_K}\lambda-K(t,z)] dt$$

Since $$\frac{d}{ds} (\psi^s)^*(\lambda)=(\psi^s)^*(L_{X_K}\lambda))= (\psi^s)^*(di_{X_K} \lambda + i_{X_K}d\lambda)= d [(\psi^s)^*(i_{X_K} \lambda+K)] ds
$$ we have

$$(\psi^t)^*\lambda - \lambda= d \int_0^t [(\psi^s)^*(i_{X_K} \lambda+K)] ds = dF(t,x)$$

and thus
$$\Psi^*(\lambda+hdt)=dF(t,z)+\frac{\partial}{\partial t}F(t,z)dt = DF(t,z)$$

Note that if $K$ is $1$ periodic in time, and $\psi^1=\psi^0=\Id$ (this implies $\psi^{t+1}=\psi^t$)we get a Lagrangian submanifold in $T^*(S^1\times N)$. This Lagrangian will be  exact, provided we change $K$ by some constant, by the following arguments:
 \begin{enumerate} 
\item $ F(t+1,z)-F(t,z)$ is constant in time, since $$\frac{d}{dt} (F(t+1,z)-F(t,z))= (\psi^{t+1})^*(i_{X_K} \lambda+K)-(\psi^{t})^*(i_{X_K}\lambda +K)$$

\item   $F(1,z)-F(0,z)$ is constant in $z$, since  $$dF(t+1,z)-dF(t,z)= (\psi^{t+1})^*(\lambda)- (\psi^{t})^*(\lambda)$$

\item  According to (a) and (b), $F(t+1,z)-F(t,z)$ is constant $c$. Since changing $K$ by a constant $c$, changes $F(1,z)-F(0,z)$ by $c$. The proof is now clear.
\end{enumerate} 
 \end{proof} 

%

%

Now let $H(t,z)$ be a Hamiltonian on $(M,\omega)$. We associate to it the Lagrangian manifold \begin{equation}\tag{1} \Lambda_H=\{(t,-H(t,\phi^t(x)), x ,\phi^t(x))\}\end{equation}

By a simple computation as above, it is indeed Lagrangian. If $\phi^t$ is $C^0$ close to the identity, by Weinstein neighbourhood's theorem $\Lambda_H$ will be contained in $T^*([0,1]\times \Delta_M)$ where $\Delta_M$ is the diagonal in $M\times \overline{M}$. As a submanifold of $T^*([0,1]\times \Delta_M)$ it will then be exact, since it is constructed as the  image of the above map $\Psi$ associated to $K(t,x_1,x_2)=H(t,x_2)$.

Moreover, according to the lemma, if near $t=0$ and $t=1$ we have both that $H$ vanishes  and that $\phi^t(z)=z$,  $\Lambda$ may be closed to a  Lagrangian submanifold of $T^*(S^1)\times T^*N$ and , after shifting $H$ by a constant, this Lagrangian will also be exact. 

This being said, to a  $C^1$ map $\chi :[0,1] \to [0,1]$ and a  Hamiltonian $H(t,z)$ with flow $\phi^t$  we associate the flow $$\phi_\chi^t(z)= \phi^{\chi(t)}(z)$$  generated by $$H_\chi(t,z)=\chi'(t)H(\chi(t),z)$$

Now if $H_n(t,z)$ is a sequence converging to $H(t,z)$ such that $\phi_n^t$ converges to the  identity map, the sequence $H_{n,\chi}(t,z)$ converges to $H_\chi(t,z)$, and $\phi_{n,\chi}^t$ converges to identity. 

We shall assume $\chi $ is identically zero in a neighbourhood of $0$ and $1$.

Consider now the Lagrangians $$\Lambda_n=\{t,-H_{n,\chi}(t,\phi_{n,\chi}^t(z)),z,\phi_{n,\chi}^t(z)) \mid t \in [0,1], z \in M\}$$
Note that since $H_\chi(0,z)=0$, these are Lagrangians of the type \thetag{1}. Since 

 $$H_{n,\chi}(t,\phi_{n,\chi}^t(z))=\chi'(t)H_n(\chi(t),z)=0$$
and
 $$\phi_{n,\chi}^t(z))=\phi_n^{\chi(t)}(z)=z$$
we may close $\Lambda_n$ to an exact Lagrangian submanifold   in $T^*(S^1)\times M \times \overline{M}$. 

The $\Lambda_n$ converge in the  $C^0$ topology to 

$$\Lambda = \{t,-H_{\chi}(t,z),z,z) \mid t \in [0,1], z \in M\} \subset T^*(S^1)\times M \times \overline{M}$$

Since $\Lambda$ is contained in $T^*(S^1) \times \Delta_M$, where $\Delta_M$ is the diagonal in$M\times \overline{M}$ we may assume using Weinstein's theorem, that $\Lambda_n$ is contained in a set symplectomorphic to a neighbourhood of the zero section in  $T^*(S^1)\times \Delta_M$. 



 
  Since the $\Lambda_n$ are exact Lagrangians, $\Lambda$ must be  Lagrangian according to Corollary 2.2. Thus there is a constant $c_\chi$ such that the form $(H_\chi(t,z)-c_\chi) dt$ must be closed in the sense of distributions. This implies that $H_\chi(t,z)=h_\chi(t)$ but since for each $t$ the average of $H_\chi$ over $M$ is zero, we must have $$\chi'(t)H(\chi(t),z)=H_\chi(t,z)=c_\chi$$ for all $\chi$ satisfying the above assumption. Since $H_\chi$ vanishes at $t=0$, we must have $c_\chi=0$. 
  Now it is not hard, for any $t_0\in ]0,1[$ to find such a $\chi$ with $\chi(t_0)=t_0$ and $\chi'(t_0)\neq 0$. This implies that $H(t_0,z)=0$. Since this holds on $]0,1[$, and $H$ is continuous, it must hold everywhere. 
  
\begin{remark}
One would like to know whether proposition 2.1 still holds for $N$ a
general compact manifold. This does not seem to follow literally from \cite{Laudenbach-Sikorav-IMRN}, even though their method may be useful. \end{remark}
 \begin{remark} 
 We could have also used the ideas from \cite{Sik3} for most of our proof. We think however that
 proposition 2.1 is of independent interest.  
 \end{remark}

\end{document}